\documentclass[11pt]{article}
\usepackage{cite}
\usepackage{amsmath,amsthm,amssymb,amsfonts}
\usepackage{amsfonts,amssymb,color}
\usepackage{graphicx}
\usepackage{float}
\usepackage{curves}
\usepackage{subfigure}
\usepackage{fullpage}

\newtheorem{theorem}{Theorem}[section]

\newtheorem{lemma}[theorem]{Lemma}
\newtheorem{corollary}[theorem]{Corollary}

\numberwithin{equation}{section}

\title{A tight lower bound on the matching number of graphs via Laplacian eigenvalues}

\author{Xiaofeng Gu\thanks{
Department of Computing and Mathematics, University of West Georgia, Carrollton, GA 30118, USA. \newline Email: {\tt xgu@westga.edu}
}\ \
and
Muhuo Liu\thanks{
Department of Mathematics, South China Agricultural University, Guangzhou, 510642, China.} 
\thanks{Research Center for Green Development of Agriculture, South China Agricultural University, Guangzhou, 510642, China.
Email: {\tt liumuhuo@163.com}}
}

\begin{document}
\date{}
\maketitle
\noindent

\begin{abstract}
Let $\alpha'$ and $\mu_i$ denote the matching number of a non-empty simple graph $G$ with $n$ vertices and the $i$-th smallest eigenvalue of its Laplacian matrix, respectively. In this paper, we prove 
a tight lower bound $$\alpha' \ge \min\left\{\Big\lceil\frac{\mu_2}{\mu_n} (n -1)\Big\rceil,\ \ \Big\lceil\frac{1}{2}(n-1)\Big\rceil \right\}.$$
This bound strengthens the result of Brouwer and Haemers who proved that if $n$ is even and $2\mu_2 \ge \mu_n$, then $G$ has a perfect matching.
A graph $G$ is factor-critical if for every vertex $v\in V(G)$, $G-v$ has a perfect matching. We also prove an analogue to the result of Brouwer and Haemers mentioned above by showing that if $n$ is odd and $2\mu_2 \ge \mu_n$, then $G$ is factor-critical. We use the separation inequality of Haemers to get a useful lemma, which is the key idea in the proofs. This lemma is of its own interest and has other applications. In particular, we prove similar results for the number of balloons, spanning even subgraphs, as well as spanning trees with bounded degree.
\end{abstract}



\section{Introduction}
Throughout this paper, we consider simple undirected graphs. A graph is non-empty if it has at least one edge.
A {\bf matching} in a graph is a set of pairwise non-adjacent edges.
A {\bf maximum matching} is one with maximum size among all matchings of the graph.
The {\bf matching number} of a graph $G$, denoted by $\alpha'(G)$, is the size of a maximum matching in $G$.
A {\bf perfect matching} in a graph $G$ with $n$ even is a matching of size $n/2$.
It has been shown by Tutte~\cite{Tutte47} that a graph $G$ has a perfect matching if and only if $o(G-S)\le |S|$ for every
vertex subset $S\subseteq V(G)$, where $o(G-S)$ denotes the number of odd components of $G-S$.
Generalizing Tutte's condition, Berge discovered the following formula for the matching number.
The formula can also be derived from Tutte's theorem, and thus is also called the {\bf Berge-Tutte Formula}.
\begin{theorem}[Berge~\cite{Berge58}]
\label{btfo}
The matching number of a graph $G$ with $n$ vertices  is $$\alpha'(G) = \frac{1}{2} \left( n+ \min_{S\subseteq V(G)} \{|S| - o(G-S)\} \right).$$
\end{theorem}

The lower bounds on the matching number of regular graphs have been studied in, among other references, \cite{BDDFK04, HeYe07, OWe10, OWe11}.
Biedl et.~al~\cite{BDDFK04} proved a tight lower bound on $\alpha'(G)$ of a cubic graph $G$. It was generalized to regular graphs
with higher degree by Henning and Yeo~\cite{HeYe07} and by O and West~\cite{OWe10, OWe11}.
The matching number of subcubic graphs (graphs with maximum degree at most 3) was well studied in \cite{BDDFK04, HeLR12, HaSc17}.
In particular, Haxell and Scott~\cite{HaSc17} gave a complete description of triples $(\alpha, \beta, \gamma)$ of real numbers for which there exists
a constant $K$ such that $\alpha'(G)\ge \alpha n_3 + \beta n_2 +\gamma n_1 -K$ holds for every connected subcubic graph $G$,
where $n_i$ denotes the number of vertices of degree $i$ for $i=1,2,3$.
Lower bounds involving other parameters have been discovered in \cite{HeLR12, OWe10, OWe11}.

The bounds of the matching
number in general case are less studied.
Henning and Yeo~\cite{HeYe18} established a tight lower bound of the matching number of a graph with given maximum degree, through which, they were able to give a complete description of pairs $(\gamma, \beta)$ of real numbers for which there exists a constant $K$ such that $\alpha'(G)\ge \gamma n + \beta m -K$ holds for every connected graph $G$ with maximum degree at most $k$, where $n$ and $m$ denote the number of vertices and the number of edges, respectively.

\par\medskip

The matrix $L=D-A$ is the {\bf Laplacian matrix} (also called {\bf combinatorial Laplacian} or {\bf discrete Laplacian}) of a graph $G$, where $D$ is the diagonal matrix of vertex degrees  and $A$ is the adjacency matrix of $G$. The eigenvalues of the Laplacian (resp. adjacency) matrix of $G$ are called
the Laplacian (resp. adjacency) eigenvalues of $G$. Throughout the paper, $\mu_i(G)$ denotes the $i$-th smallest Laplacian  eigenvalue of   $G$. 
It holds that $\mu_1(G)=0$. The second smallest Laplacian eigenvalue $\mu_2(G)$  is known as the {\bf algebraic connectivity} of $G$.
We have $\mu_2(G)=0$ if and only if $G$ is disconnected. In what follows, if there is no risk of confusion, we always simplify $\mu_i(G)$ as $\mu_i$.

The studies of perfect matchings and the matching number of regular graphs by means of the adjacency eigenvalues have appeared in the literature \cite{BrHa05, Cioa05, CiGr07, CiGH09, KrSu06, OCi10}, among others. Brouwer and Haemers~\cite{BrHa05} initiated the research of sufficient conditions for perfect matchings of regular graphs in terms of the adjacency eigenvalues, followed by Cioab\u{a} and Gregory~\cite{CiGr07} who extended the idea to the matching number of regular graphs. Cioab\u{a}, Gregory and Haemers~\cite{CiGH09} further provided a tight sufficient condition involving the third largest adjacency eigenvalue of a connected regular graph, which was then extended to regular graphs with fixed edge-connectivity by O and Cioab\u{a}~\cite{OCi10}.

For graphs that are not necessarily regular, some results have been obtained in \cite{FYZ07, O21} via the largest adjacency eigenvalues. However, the methods via other adjacency eigenvalues used in the literature~\cite{BrHa05, Cioa05, CiGr07, CiGH09, KrSu06, OCi10} rely much on the regularity of graphs, and most results seem not possible to be extended to irregular cases. By using Laplacian eigenvalues, Brouwer and Haemers~\cite{BrHa05} discovered a sufficient condition for perfect matchings. In the following, we always assume that $G$ is a graph with $n$ vertices and at least one edge.

\begin{theorem}[Brouwer and Haemers~\cite{BrHa05}]\label{newt}
If $n$ is even and   $2\mu_2 \ge \mu_n$, then $G$ has a perfect matching.
\end{theorem}

In this paper,  we generalize Theorem~\ref{newt} from perfect matchings to the matching number, and discover a tight lower bound on the matching number of graphs in terms of the Laplacian eigenvalues. 

\begin{theorem}
\label{amnu}
Let $r$ be a real number with $0< r\le 1/2$. 
If $\mu_2 \ge r \mu_n,$ then $\alpha'(G)\ge r(n-1)$.
\end{theorem}
The condition ``$\mu_2 \ge r \mu_n$'' of Theorem~\ref{amnu} is tight. To see this, let $G=K_{s,t}$ be a complete bipartite graph with $n=s+t$ vertices, where $t = \lceil \frac{1-r}{r}s +\frac{a}{r} \rceil$ for some small fixed real number $a$ with  $a > r$. Then $t >s$, $\mu_2 =s$ and $\mu_n =n= s+t  \ge \frac{s+a}{r}$. Thus, $\alpha'(K_{s,t}) =s \le rn -a < r(n-1)$. However, $\mu_2/\mu_n = s/(s+t)\leq rs/(s+a) < r$, and by choosing $s$ to be sufficiently large, we can see that $\mu_2/\mu_n$ is arbitrarily close to $r$.

The condition in Theorem~\ref{amnu} is also tight for non-bipartite graphs. 
For a positive integer $s$, let $t = \lceil \frac{1-r}{r}s +\frac{a}{r}\rceil > s$ for some small fixed real number $a$ with $a > r$.
Let $G=K_s\vee tK_1$ be the join of the complete graph $K_s$ on $s$ vertices and  the union of $t$ isolated vertices $tK_1$.
That is, $K_s\vee tK_1$ is the graph obtained by joining every vertex of $K_s$ to every vertex of $tK_1$.
Then the number of vertices of $G$  is $n = s+t \ge \frac{s+a}{r}$, $\mu_2 =s$ and  $\mu_n =n= s+t$.
Thus $\alpha'(K_s\vee tK_1) =s \le rn -a < r(n -1)$.
However, $\mu_2/\mu_n = s/(s+t)\leq rs/(s+a) < r$ can be arbitrarily close to $r$.

\par\medskip 
By choosing $ r=\min\{\frac{\mu_2}{\mu_n}, \frac{1}{2}\}$, 
we have the following tight lower bound on $\alpha'(G)$.

\begin{corollary}
$\displaystyle \alpha'(G)\ge \min\left\{\Big\lceil\frac{\mu_2}{\mu_n} (n -1)\Big\rceil,\ \ \Big\lceil\frac{1}{2}(n-1)\Big\rceil \right\}.$
\end{corollary}

A graph $G$ is {\bf factor-critical} if for every vertex $v\in V(G)$, $G-v$ has a perfect matching. Similar to Tutte's perfect matching theorem,
Gallai~\cite{Gall63} proved that $G$ is factor-critical if and only if $|V(G)|$ is odd and $o(G-S)\le |S|$ for every non-empty subset $S\subseteq V(G)$.
Spectral conditions of factor-critical regular graphs have been studied in \cite{CiGr07, CiGH09}.
By using Gallai's theorem, we can get the following tight sufficient condition for factor-critical graphs. 
\begin{theorem}
\label{criti}
If $n$ is odd and $2\mu_2 \ge \mu_n$, then $G$ is factor-critical.
\end{theorem}
To see that the result of Theorem~\ref{criti} is sharp, observe  that $K_{s, s+1}$ (in fact, no bipartite graph is factor-critical) and $K_s\vee (s+1)K_1$
are not factor-critical, but $2\mu_2/\mu_n = 2s/(2s+1)<1$ and can be arbitrarily close to $1$ by choosing sufficiently large $s$.

\par\medskip
In the next section, we will prove a useful lemma as well as Theorems~\ref{amnu} and \ref{criti}. As applications of this new lemma, in Sections~\ref{snb} and \ref{sss}, we will show similar results for the number of balloons, spanning even subgraphs, and spanning trees with degree bounded above by a fixed number.

\section{Proofs of a useful lemma and main theorems}

Our main tool to connect graph parameters to eigenvalues is the following theorem.
\begin{theorem}[Haemers~\cite{Haem95}]
Let $G$ be a graph with $n$ vertices and at least one edge. Suppose that $X$ and $Y$ are two  disjoint subsets of $V(G)$ such that there is no edge between $X$ and $Y$. Then
\begin{equation}
\label{muxy}
\frac{|X| |Y|}{(n-|X|)(n-|Y|)} \le \left(\frac{\mu_n - \mu_2}{\mu_n + \mu_2}\right)^2.
\end{equation}
\end{theorem}

\par\medskip
Utilizing the above separation inequality, we can prove a simple yet useful lemma, which is essential to our proofs.
The lemma is of its own interest and we will present more applications of this lemma in the next two sections.
\begin{lemma}
\label{22l}
Let $G$ be a graph with $n$ vertices and at least one edge. Suppose that $S\subset V(G)$ such that $G-S$ is disconnected.
Let $X$ and $Y$ be disjoint vertex subsets of $G-S$ such that $X\cup Y= V(G)-S$ with $|X|\le |Y|$. Then
\begin{equation}
\label{xupp}
|X| \le \frac{\mu_n - \mu_2}{2\mu_n}\cdot n,
\end{equation}
and
\begin{equation}
\label{ssiz}
|S| \ge \frac{2\mu_2}{\mu_n -\mu_2}\cdot |X|,
\end{equation}
with each equality holding  only when $|X|=|Y|$.
\end{lemma}
\begin{proof}
By (\ref{muxy}), we have
\begin{equation}
\label{prodxy}
|X| |Y| \le \left(\frac{\mu_n - \mu_2}{\mu_n + \mu_2}\right)^2 (n-|X|)(n-|Y|),
\end{equation}
Let $\beta = \frac{\mu_n - \mu_2}{\mu_n + \mu_2}$. Then $0<\beta \le 1$, as $G$ contains at least one edge and $G-S$ is disconnected. It follows that
$$ |X|^2 \le |X|\cdot |Y| \le \beta^2 (n-|X|)(n-|Y|)\le \beta^2 (n-|X|)^2,$$
that is $$|X| \le \beta (n - |X|),$$ and hence
\begin{equation}
\label{betax}
|X|\le \frac{\beta n}{1+\beta}=\frac{\mu_n - \mu_2}{2\mu_n}\cdot n,
\end{equation}
with the equality holding  only when $|X|=|Y|$.

\par\medskip
Also, since $|Y| = n - |S| - |X|$, by (\ref{prodxy}), we have
$$ |X| (n - |S| - |X|) = |X|\cdot |Y| \le \beta^2 (n - |X|)(n - |Y|) = \beta^2 (n - |X|)(|S| + |X|),$$
implying that
\begin{equation}
\label{xnbo}
|X| n \le \left(\beta^2 (n - |X|) + |X| \right) \left(|S| + |X| \right)
= \left(\beta^2 n + (1-\beta^2)|X| \right) \left(|S| + |X| \right).
\end{equation}
By (\ref{betax}), we have
$$(1-\beta^2)|X|\le (1-\beta^2) \cdot \frac{\beta n}{1+\beta} = (\beta - \beta^2) n,$$
which, together with (\ref{xnbo}), implies that
$$|X|n \le \left(\beta^2 n + (\beta - \beta^2) n \right) \left(|S| + |X| \right) = \beta n \left(|S| + |X| \right),$$
and we have
$$ |X| \le \beta \left(|S| + |X| \right).$$
Hence,
\begin{equation*}
|S| \ge \frac{1 -\beta}{\beta} |X| = \frac{2\mu_2}{\mu_n - \mu_2}\cdot |X|.
\end{equation*}
Since (\ref{betax}) was utilized, the equality holds in (\ref{ssiz}) only when $|X|=|Y|$.
\end{proof}

\par\medskip
\begin{proof}[\bf Proof of Theorem~\ref{amnu}]
By Theorem~\ref{btfo}, it suffices to show that for every $S\subseteq V(G)$, $$|S| - o(G-S)\ge (2r-1)n -2r,$$
i.e., $$o(G-S) - |S|\le (1-2r)n +2r.$$

If $o(G-S)\le 1$, then $o(G-S) - |S| \le 1- |S| \le 1 \le (1-2r)n +2r$, since $1-2r\le (1-2r)n$. Thus we may assume that $o(G-S)\ge 2$.
Let $O_1, O_2, \ldots, O_c$ be vertex sets of the odd components of $G-S$, where $c=o(G-S)$.
Without loss of generality, suppose that $|O_1|\le |O_2| \le \cdots \le |O_c|$.

We first show that if $|V(G) -S| = c$, then the theorem holds. In this case, each $O_i$ is an isolated vertex for $1\le i\le c$.
Then $|S|= n-c$. Define $X=\bigcup_{1\leq i\leq \lfloor c/2\rfloor}O_i$ and $Y=V(G)-S-X$. By (\ref{ssiz}),
\begin{equation*}
n-c= |S| \ge \frac{2\mu_2}{\mu_n -\mu_2}\cdot |X|\geq \frac{2\mu_2}{\mu_n -\mu_2}\cdot \frac{c-1}{2}
=\frac{(c-1)\mu_2}{\mu_n -\mu_2},
\end{equation*}
which implies that
\begin{equation}\label{ne}
c \le \frac{\mu_n - \mu_2}{\mu_n} n + \frac{\mu_2}{\mu_n}.
\end{equation}
Combining (\ref{ne}) with  $\mu_2 \ge r\mu_n$,
\begin{equation*}
o(G-S)-|S| = 2c-n \le 2\left(\frac{\mu_n - \mu_2}{\mu_n} n + \frac{\mu_2}{\mu_n}\right) - n
= n -\frac{2\mu_2}{\mu_n}(n-1) \le n - 2r(n-1) = (1-2r)n +2r,
\end{equation*}
as desired.

Thus, we may assume that $|V(G)-S| \ge c+1$ in the following.
We claim that $O_1, O_2, \ldots, O_c$ can be partitioned into two sets $X$ and $Y$ such that $|Y|\ge |X|\ge c/2$.
If $c$ is even, we can simply define $X =\bigcup_{1\le i\le \lfloor c/2\rfloor} O_i$ and $Y = V(G) - S -X$.
Now we assume that $c$ is odd. If $|O_{(c-1)/2}|\ge 2$, then define $X =\bigcup_{1\le i\le (c-1)/2} O_i$ and $Y = V(G) - S -X$ as needed.
The remaining case is $|O_1| = \cdots = |O_{(c-1)/2}|=1$. We can define $X =\bigcup_{1\le i\le (c+1)/2} O_i$ and $Y = V(G) - S -X$,
and we need to show that $|Y|\ge |X|\ge c/2$. If $|O_{(c+1)/2}| =1$,
then $|X| = \frac{c+1}{2}$ and $|Y| = |V(G)-S| -|X| \ge \frac{c+1}{2}$, since $|V(G)-S| \ge c+1$.
If $|O_{(c+1)/2}| \ge 2$, then $|X| = \frac{c-1}{2} + |O_{(c+1)/2}| \ge \frac{c-1}{2} +2 > \frac{c}{2}$
and $|Y| = \sum_{i > (c+1)/2} |O_i| \ge 2\cdot \frac{c-1}{2}=c-1\ge c/2$. Switch $X$ and $Y$ whenever needed to get $|Y|\ge |X|$.

\par\medskip
It follows that $o(G-S) = c\le 2|X|$. Thus, by (\ref{ssiz}),
\begin{equation}
\label{cceven}
o(G-S) - |S| = c -|S| \le 2|X| -|S| \le 2|X| - \frac{2\mu_2}{\mu_n -\mu_2} |X| =  \frac{2(\mu_n - 2\mu_2)}{\mu_n -\mu_2}|X|.
\end{equation}

If $2\mu_2 > \mu_n$, then $o(G-S) -|S| <0 < (1-2r)n +2r$, done.
Thus we may assume that $2\mu_2 \le \mu_n$. Since $\mu_2 \ge r\mu_n$, we have, by (\ref{xupp}) and (\ref{cceven}),
$$o(G-S) -|S| \le  \frac{2(\mu_n -2\mu_2)}{\mu_n -\mu_2} \cdot \frac{\mu_n - \mu_2}{2\mu_n}\cdot n = \frac{\mu_n - 2\mu_2}{\mu_n}\cdot n \le (1-2r)n,$$
and so $o(G-S) - |S| < (1-2r)n +2r$, as desired. 
\end{proof}

\par\medskip
\begin{proof}[\bf Proof of Theorem~\ref{criti}]
Gallai~\cite{Gall63} proved that $G$ is factor-critical if and only if $|V(G)|$ is odd and $o(G-S)\le |S|$ for every non-empty subset $S\subseteq V(G)$.
Suppose by way of contradiction that $G$ is not factor-critical. Then there exists a non-empty subset $S\subseteq V(G)$ such that
\begin{equation}
\label{ovvs}
o(G-S)\ge |S| +1 \ge 2.
\end{equation}

Let $O_1, O_2, \ldots, O_c$ be vertex sets of the odd components of $G-S$, where $c=o(G-S)$.
Without loss of generality, suppose that $|O_1|\le |O_2| \le \cdots \le |O_c|$.
Define $X =\bigcup_{1\le i\le \lfloor c/2\rfloor} O_i$ and $Y = V(G) - S -X$. Then $\frac{c-1}{2} \le |X|\le |Y|$.
By (\ref{ssiz}) and $|X|\geq \frac{c-1}{2}$, we have
\begin{equation}
\label{sscc}
|S| \ge  \frac{2\mu_2}{\mu_n - \mu_2}\cdot (c-1)/2,
\end{equation}
with equality holding only when $\frac{c-1}{2}=|X|=|Y|$. However, if $c$ is even, then $\frac{c-1}{2}<\frac{c}{2}\le |X|$; and
if $c$ is odd, then $|X| < |Y|$ by definitions of $X$ and $Y$. Thus the equality never hold in (\ref{sscc}). Since $2\mu_2 \ge \mu_n$, from (\ref{sscc}) we can deduce that
\begin{equation}
|S| >  \frac{2\mu_2}{\mu_n - \mu_2}\cdot (c-1)/2 \ge c-1 = o(G-S) -1,
\end{equation}
contrary to (\ref{ovvs}). This completes the proof.
\end{proof}

\section{The number of balloons in a graph}
\label{snb}
Defined in \cite{OWe10}, a {\bf balloon} in a graph $G$ is a maximal 2-edge-connected subgraph of $G$ incident to exactly one cut-edge of $G$.
In other words, a balloon is joined to the rest of the graph by a single edge. Note that a balloon may contain cut-vertices and blocks. Clearly maximal 2-edge-connected subgraphs in a graph are pairwise vertex-disjoint. The number of balloons in $G$ is denoted by $b(G)$.  This parameter has been shown to be closely related to the matching number in \cite{OWe10}. Actually, by proving an upper bound on the number of balloons in a connected odd regular graph $G$, O and West~\cite{OWe10} obtained a lower bound on $\alpha'(G)$. The upper bound on $b(G)$ of a regular graph $G$ has been studied in \cite{OCi10}  by means of the adjacency eigenvalues. We aim to discover an upper bound on $b(G)$ of a graph $G$ that is not necessarily regular in terms of its Laplacian eigenvalues.

\begin{figure}[htb]
\centering
\subfigure[A graph with $b(G)=3$]{\includegraphics[width=0.33\textwidth]{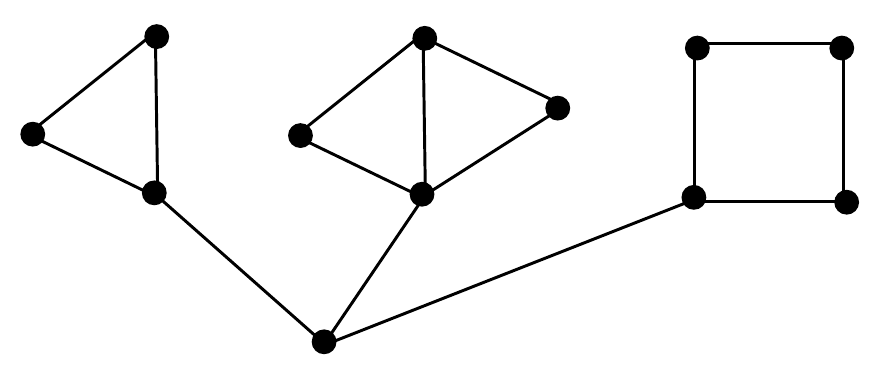}}
\hspace{.8in}
\subfigure[A dumbbell]{\includegraphics[width=0.28\textwidth]{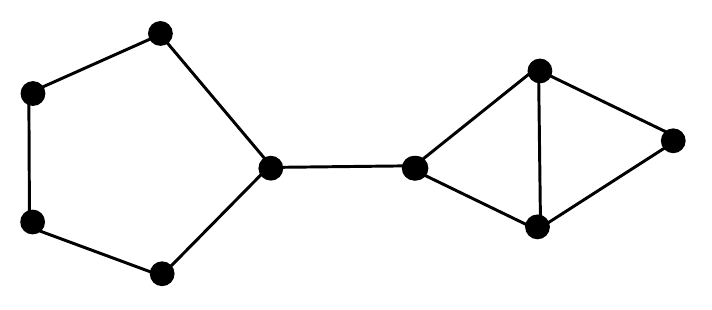}}
\caption{Examples of balloons}
\label{fig:balloon}
\end{figure}

Let $G$ be a connected  graph with minimum degree $\delta$. Then every balloon of $G$ must have at least $\max\{3, \delta +1\}$ vertices, since it is 2-edge connected and contains a vertex of degree at least $\max\{2, \delta\}$. Thus $b(G)\le \min\{n/3, n/(\delta +1)\}$ trivially.
A {\bf dumbbell} is a graph consisting of exactly two balloons joined by a single edge. Two examples are shown in Figure~\ref{fig:balloon}.
If $G$ is a dumbbell, then $b(G)=2$. Thus it suffices to consider non-dumbbell graphs.

\begin{theorem}
Let $G$ be a  connected  non-dumbbell graph on $n$ vertices with minimum degree $\delta$,
and let $r$ be a real number with $0< r\le \min\{\frac{1}{3}, \frac{1}{\delta +1}\}$.
\\$(i)$ If $\delta \ge 2$ and $\mu_2 \ge (1 - r -r\delta) \mu_n,$ then $b(G)\le \lceil rn\rceil$.
\\$(ii)$ If $\delta =1$ and $\mu_2 \ge (1-3r) \mu_n,$ then $b(G)\le  \lceil rn\rceil$.
\end{theorem}
\begin{proof}
(i) Suppose the contrary that $b(G) \ge  \lceil rn\rceil +1 \ge 2$. Notice that all balloons are pairwise vertex-disjoint.
Since $G$ is not a dumbbell, no edge joins a pair of balloons of $G$.
Let $V_1, V_2, \ldots, V_b$ be the  vertex sets of  the balloons of $G$, and let $S$ be the set of the remaining vertices of $G$.
Without loss of generality, suppose that $|V_1|\le |V_2| \le \cdots \le |V_b|$.
Define $X =\bigcup_{1\le i\le \lfloor b/2\rfloor} V_i$ and $Y = V(G) - S -X$. Then $\frac{b-1}{2}\cdot (\delta +1) \le |X|\le |Y|$, since
each balloon contains at least $\delta +1$ vertices.
Now, by $b(G) \ge  rn+1$, $\mu_2 \ge (1 - r -r\delta)\mu_n$ and (\ref{xupp}), we have
\begin{equation}\label{n3e}
\frac{r(\delta +1)n}{2}\le \frac{b-1}{2}\cdot (\delta +1) \le |X| \le \frac{\mu_n - \mu_2}{2\mu_n}\cdot n \le \frac{r(\delta +1)n}{2},
\end{equation}
and thus all equalities hold in (\ref{n3e}). This implies that $\frac{b-1}{2}\cdot (\delta +1) =|X|=|Y|$, since (\ref{xupp}) was employed.
However, if $b$ is even, then $\frac{b-1}{2}\cdot (\delta +1) <\frac{b}{2}\cdot (\delta +1)  \le |X|$; and
if $b$ is odd, then $|X| < |Y|$ by definitions of $X$ and $Y$. Thus the equality never hold in (\ref{n3e}), a contradiction.
Hence $b(G)\le \lceil rn\rceil$.

(ii) When $\delta =1$, the proof is quite similar to part (i), except that in this case, each balloon contains at least 3 vertices,
instead of $\delta+1$ vertices. Thus the proof is omitted.
\end{proof}

\section{Spanning subgraphs with degree constraints}
\label{sss}
A perfect matching can be seen as a spanning subgraph with degree of each vertex is exactly one. Since not every graph has a perfect matching, we are interested in determining a spanning subgraph with bounded degree, i.e., the degree of each vertex is bounded from below and/or from above by some fixed number. Spanning subgraphs in which every vertex has degree $k$ of regular graphs have been studied by means of the adjacency eigenvalues in \cite{Gu14, Lu10, Lu12}. However, there are no similar results for irregular graphs. 
We can extend the results on matchings  to spanning subgraphs with bounded degree in terms of the Laplacian eigenvalues. In fact, we can do better by requiring the  parity of degree or connectedness of subgraphs. A graph $G$ is {\bf even} if the degree of each vertex is positive even. The study of even subgraphs, closely related to Eulerian graphs, has been well developed and has wide applications in graph theory.

\begin{theorem}
Let $G$ be a graph on $n$ vertices with minimum degree $\delta\ge 3$.
If $(\delta-1)\mu_2 \ge \mu_n$, then $G$ has an even spanning subgraph.
\end{theorem}
\begin{proof}
It is mentioned in Theorem 6.2 of \cite{AkKa11} that a graph $G$ has a spanning even subgraph if and only if for every vertex subset
$S\subseteq V(G)$, $$q(G-S)\le \sum_{v\in S} \big(d_G(v) -2\big),$$ where $q(G-S)$ is the number of components $C$ of $G-S$
such that the number of edges between $C$ and $S$ is odd and $d_G(v)$ is the degree of vertex $v$ in $G$.

Suppose by way of contradiction that $G$ has no spanning even subgraphs. Then there exists a subset $S\subseteq V(G)$ such that
\begin{equation}
\label{qvs}
q(G-S)\ge \sum_{v\in S} \big(d_G(v) -2\big) +1 \ge (\delta -2)|S| +1.
\end{equation}
Clearly $S$ is not empty, since otherwise, $q(G-S)=0$ by definition of $q(G-S)$. Thus $|S|\ge 1$ and so
$q(G-S)\ge(\delta -2)|S| +1 \ge \delta-1 \ge 2$.

Let $O_1, O_2, \ldots, O_c$ be the  vertex sets of the components $C$ of $G-S$ such that the number of edges
between $C$ and $S$ is odd, where $c=q(G-S)\ge 2$.
Without loss of generality, suppose that $|O_1|\le |O_2| \le \cdots \le |O_c|$.
Define $X =\bigcup_{1\le i\le \lfloor c/2\rfloor} O_i$ and $Y = V(G) - S -X$. Then $\frac{c-1}{2} \le |X|\le |Y|$.
By (\ref{ssiz}) and $|Y|\geq |X|\geq (c-1)/2$, we have
\begin{equation}\label{n2e}
|S| \ge \frac{2\mu_2}{\mu_n -\mu_2}\cdot |X| \ge \frac{(c-1)\mu_2}{\mu_n -\mu_2},
\end{equation}
with equality holding only when $\frac{c-1}{2}=|X|=|Y|$. However, if $c$ is even, then $\frac{c-1}{2}<\frac{c}{2}\le |X|$; and
if $c$ is odd, then $|X| < |Y|$ by definitions of $X$ and $Y$. Thus the equalities of (\ref{n2e}) never hold simultaneously.
Since $(\delta-1)\mu_2 \ge \mu_n$, (\ref{n2e}) implies that
$$|S| > \frac{(c-1)\mu_2}{\mu_n -\mu_2} \ge \frac{c-1}{\delta -2}=\frac{q(G-S) -1}{\delta -2},$$
and so $q(G-S) < (\delta -2)|S| + 1,$ contrary to (\ref{qvs}).
\end{proof}

We also have the following sufficient condition for the existence of a spanning tree with degree bounded above by a fixed number $k$,
which equivalently implies a connected spanning subgraph with degree at most $k$.
When $k=2$, such a spanning tree is exactly a Hamiltonian path, and a sufficient condition has been given by, among others,
Butler and Chung~\cite{BuCh10}, who borrowed the idea from \cite{KrSu03} (both \cite{KrSu03, BuCh10} studied a stronger
structure, i.e., Hamiltonian cycle). Thus we focus on $k\ge 3$. This work is actually motivated by \cite{CiGu16} which studied the same problem
for regular graphs and posed a question to extend the result to all simple graphs and Laplacian eigenvalues.
By using Lemma~\ref{22l}, we give an answer to this question in the following theorem.
\begin{theorem}
\label{bspt}
Let  $k\ge 3$ be an integer and $G$ be a connected graph on $n\ge 2$ vertices. If $(k-1)\mu_2 \ge \mu_n$,
then $G$ has a spanning tree with maximum degree at most $k$.
\end{theorem}
\begin{proof}
Let  $k\ge 3$. Win~\cite{Win89}   proved the following result for connected graphs: If for every
subset  $S\subseteq V(G)$, we have  $c(G-S)\le (k-2)|S| +2$,
then $G$ has a spanning tree with maximum degree at most $k$, where $c(G-S)$ is the number of components in $G-S$
(see also~\cite{ElZh00}).

Suppose   the contrary that $G$ has no such spanning trees. Then there is a vertex subset $S\subseteq V(G)$ such that
\begin{equation}
\label{cgst}
c(G-S)\ge (k-2)|S| +3.
\end{equation}

Let $V_1, V_2, \ldots, V_c$ be vertex sets of the components of $G-S$, where $c=c(G-S)\ge 3$.
Without loss of generality, suppose that $|V_1|\le |V_2| \le \cdots \le |V_c|$.
Define $X =\bigcup_{1\le i\le \lfloor c/2\rfloor} V_i$ and $Y = V(G) - S -X$. Then $\frac{c-1}{2} \le |X|\le |Y|$.
By (\ref{ssiz}) and $(k-1)\mu_2 \ge \mu_n$, we have
\begin{equation*}
|S| \ge \frac{2\mu_2}{\mu_n - \mu_2}\cdot |X| \ge \frac{(c-1)\mu_2}{\mu_n - \mu_2} \ge \frac{c-1}{k-2},
\end{equation*}
which implies that $$(k-2)|S| +1 \ge c = c(G-S),$$ contrary to (\ref{cgst}), completing the proof.
\end{proof}

\par\medskip
\noindent
{\bf Acknowledgment.}
\\The authors would like to thank the referees for their valuable suggestions and comments which lead to an improvement of the original manuscript.
X. Gu is partially supported by a grant from the Simons Foundation (522728).



\begin{thebibliography}{99}




\bibitem{AkKa11}
J. Akiyama and M. Kano, {\it Factors and Factorizations of Graphs, Proof Techniques in Factor Theory}, springer, 2011.



\bibitem{Berge58}
C. Berge, Sur le couplage maximum d'un graphe, C.R. Acad. Sci. Paris 247 (1958) 258--259.

\bibitem{BDDFK04}
T. Biedl, E.D. Demaine, C.A. Duncan, R. Fleischer and S.G. Kobourov, Tight bounds on maximal and maximum matchings,
Discrete Math. 285 (2004), 7--15.



\bibitem{BrHa05}
A.E. Brouwer and W.H. Haemers, Eigenvalues and perfect matchings, Linear Algebra Appl. 395 (2005), 155--162.


\bibitem{BuCh10}
S. Butler and F. Chung, Small spectral gap in the combinatorial Laplacian implies Hamiltonian,
Ann. Comb. 13 (2010), 403--412.

\bibitem{Cioa05}
S.M. Cioab\u{a}, Perfect matchings, eigenvalues and expansion, C. R. Math. Rep. Acad. Sci. Canada,  27 (2005), 101--104.

\bibitem{CiGr07}
S.M. Cioab\u{a} and D.A. Gregory, Large matchings from eigenvalues, Linear Algebra Appl. 422 (2007), 308--317.

\bibitem{CiGH09}
S.M. Cioab\u{a}, D.A. Gregory and W.H. Haemers, Matchings in regular graphs from eigenvalues,
J. Combin. Theory Ser. B, 99 (2009), 287--297.

\bibitem{CiGu16}
S.M. Cioab\u{a} and X. Gu, Connectivity, toughness, spanning trees of bounded degrees, and spectrum of regular graphs,
Czechoslovak Math. J. 66 (2016), 913--924.




\bibitem{ElZh00}
M.N. Ellingham and X. Zha, Toughness, trees, and walks, J. Graph Theory 33 (2000), 125--137.

\bibitem{FYZ07}
L. Feng, G. Yu, and X.-D. Zhang, Spectral radius of graphs with given matching number, Linear Algebra Appl. 422 (2007), 133--138.

\bibitem{Gall63}
T. Gallai, Neuer Beweis eines Tutte'schen Satzes, Magyar Tud. Akad. Mat. Kut. Int. K\"ozl. 8 (1963), 135--139.

\bibitem{Gu14}
X. Gu, Regular factors and eigenvalues of regular graphs, European J. Combin. 42 (2014), 15--25.

\bibitem{Haem95}
W.H. Haemers, Interlacing eigenvalues and graphs, Linear Algebra Appl. 226-228 (1995), 593--616.

\bibitem{HaSc17}
P.E. Haxell and A.D. Scott, On lower bounds for the matching number of subcubic graphs, J. Graph Theory 85 (2017), 336--348.

\bibitem{HeLR12}
M.A. Henning, C. L\"owenstein, and D. Rautenbach, Independent sets and matchings in subcubic graphs, Discrete Math. 312 (2012), 1900--1910.

\bibitem{HeYe07}
M.A. Henning and A. Yeo, Tight lower bounds on the size of a maximum matching in a regular graph, Graphs Combin. 23 (2007), 647--657.

\bibitem{HeYe18}
M.A. Henning and A. Yeo, Tight lower bounds on the matching number in a graph with given maximum degree, J. Graph Theory 89 (2018), 115--149.


\bibitem{KrSu03} M. Krivelevich and B. Sudakov, Sparse pseudo-random graphs are Hamiltonian, J. Graph Theory 42 (2003), 17--33.

\bibitem{KrSu06}
M. Krivelevich and B. Sudakov, Pseudo-random graphs, More Sets, Graphs and Numbers,
Bolyai Soc. Math. Stud., 15, Springer, Berlin, 2006, 199--262.

\bibitem{Lu10}
H. Lu, Regular factors of regular graphs from eigenvalues, Electron. J. Combin. 17 (2010), \#R159.

\bibitem{Lu12}
H. Lu, Regular graphs, eigenvalues and regular factors, J. Graph Theory 69 (2012), 349--355.

\bibitem{O21}
S. O, Spectral radius and matchings in graphs, Linear Algebra Appl. 614 (2021), 316--324.

\bibitem{OCi10}
S. O and S.M. Cioab\u{a}, Edge-connectivity, eigenvalues and matchings in regular graphs, SIAM J. Discrete Math. 24 (2010), 1470--1481.


\bibitem{OWe10}
S. O and D.B. West, Balloons, cut-edges, matchings, and total domination in regular graphs of odd degree, J. Graph Theory 64 (2010), 116--131.

\bibitem{OWe11}
S. O and D.B. West, Matching and edge-connectivity in regular graphs, European J. Combin. 32 (2011), 324--329.


\bibitem{Tutte47}
W.T. Tutte, The factorization of linear graphs, J. London Math. Soc. 22 (1947), 107--111.


\bibitem{Win89} S. Win, On a connection between the existence of $k$-trees and the toughness of a graph, Graphs Combin. 5 (1989), 201--205.




\end{thebibliography}
\end{document}